\documentclass[twocolumn]{autart}    
\usepackage{amssymb}
\usepackage{amsmath}
\usepackage{algorithm}
\usepackage{algorithmic}
\usepackage{booktabs}
\usepackage{booktabs}
\usepackage{booktabs}
\usepackage{graphicx}  
\usepackage{subfigure}
\usepackage{epstopdf}
\usepackage{subfloat}
\usepackage{hyperref}

\newtheorem{theorem}{Theorem}
\newtheorem{corollary}{Corollary}
\newtheorem{lemma}{Lemma}
\newtheorem{remark}{Remark}
\newtheorem{assumption}{Assumption}
\newtheorem{proposition}{Proposition}
\newtheorem{example}{Example}
\newtheorem{definition}{Definition}

\newcommand{\opt}{^{\star}} 
\newcommand{\tr}{\intercal} 

\usepackage{color}
\newcommand{\change}[1]{\textcolor{black}{#1}}
\newcommand{\StartChange}{\color{black}}
\newcommand{\EndChange}{\color{black}}

\definecolor{YJ}{rgb}{0.0,0.5,0.0}

\begin{document}

\begin{frontmatter}
    
\title{Parallel Explicit Model Predictive Control\thanksref{footnoteinfo}}
    
\thanks[footnoteinfo]{This paper is extended version of the paper~\cite{Oravec2017} presented at IFAC World Congress 2017 in Toulouse. Corresponding author J.~Oravec. Tel. +421 259 325 364. Fax +421 259 325 340.}
    
\author[ShanghaiTech,SIMIT,CAS]{Yuning Jiang}\ead{jiangyn@shanghaitech.edu.cn},
\author[STU]{Juraj Oravec}\ead{juraj.oravec@stuba.sk},
\author[ShanghaiTech]{Boris Houska}\ead{borish@shanghaitech.edu.cn},
\author[STU]{Michal Kvasnica}\ead{michal.kvasnica@stuba.sk}
\address[ShanghaiTech]{School of Information Science and Technology, ShanghaiTech University, Shanghai 201210, China.}
\address[SIMIT]{Shanghai Institute of Microsystem and Information Technology, Chinese Academy of Sciences,  Shanghai 200050, China}
\address[STU]{Institute of Information Engineering, Automation, and Mathematics, Faculty of Chemical and Food Technology, Slovak University~of~Technology in Bratislava, Radlinsk\'{e}ho 9, SK812-37 Bratislava, Slovak Republic.}
\address[CAS]{University of Chinese Academy of Sciences, Beijing 100049, China}


\begin{abstract}
This paper is about a real-time model predictive control (MPC) algorithm for large-scale, structured linear \change{systems} with polytopic state and control constraints. The proposed controller receives the current state measurement as an input and computes a sub-optimal control reaction by evaluating a finite number of piecewise affine functions that correspond to the explicit solution maps of small-scale parametric quadratic programming (QP) problems. We provide \change{recursive feasibility} and \change{asymptotic stability} guarantees, which can both be verified offline. The feedback controller is suboptimal on purpose because we are enforcing real-time requirements assuming that it is impossible to solve the given large-scale QP in the given amount of time. In this context, a key contribution of this paper is that we provide a \change{bound on the sub-optimality of the controller}. Our numerical simulations illustrate that the proposed explicit real-time scheme easily scales up to systems with hundreds of states and long control horizons, system sizes that are completely out of the scope of existing, non-suboptimal Explicit MPC controllers.
\end{abstract}

\end{frontmatter}

\section{Introduction}
\label{sec:introduction}
The advances of numerical optimization methods over the last decades~\cite{Diehl2009}, in particular, the development of efficient quadratic programming problem (QP) solvers~\cite{Ferreau2014}, have enabled numerous industrial applications of \change{MPC~\cite{Qin2003}.} Modern real-time optimization and control software packages~\cite{Houska2011,Mattingley2009} achieve run-times in the milli- and microsecond range by generating efficient and reliable C-code that implements problem-tailored MPC algorithms~\cite{Diehl2002,Zavala2009}. However, as much as these algorithms and codes perform well on desktop computers or other devices with comparable computation power, the number of successful implementations of MPC on embedded industrial hardware such as programmable logic controllers (PLC) and field-programmable gate arrays (FPGA) remains limited~\cite{Ingole2015}. Here, the main question is what can be done if an embedded device has simply not enough computational power or storage space to solve the exact MPC problem online and in real-time.

Many researchers have attempted to address this question. For example, the development of Explicit MPC~\cite{oberdieck2016multi} aims at reducing both the online run-time and the memory footprint of MPC by optimizing pre-computed solution maps of multi-parametric optimization \change{problems}. However, Explicit MPC has the disadvantage that the number of polytopic regions over which the piecewise affine solution map of a parametric linear or quadratic programming problem is defined, grows, in the worst case, exponentially with the number of constraints. Some authors~\cite{GM02} have suggested addressing this issue by simplifying the MPC problem formulation, e.g., by using move-blocking~\cite{Cagienard2007}, but the associated control reactions can be sub-optimal by a large margin. Other authors~\cite{Kva:regionless:2015} have worked on reducing the memory footprint of Explicit MPC---certainly making considerable progress yet failing to meet the requirement of many practical systems with more than just a few states. In fact, despite all these developments in Explicit MPC, these methods are often applicable to problems of modest size only. As soon as one attempts to scale up to larger systems, an empirical observation is that Explicit MPC is often outperformed by iterative online solvers such as active set~\cite{Ferreau2014} or interior point methods for MPC~\cite{Mattingley2009}. \change{In this context, we mention that~\cite{KL18} has recently proposed a heuristic for reducing the number of regions of Explicit MPC by using neural-network approximations. The corresponding controller does, however, only come along with guarantees on the feasibility, stability, and performance if the exact explicit solution map happens to be recovered by the deep-learning approach.}



\change{A recent trend in optimization based control is to solve large MPC problems by breaking them into smaller ones. This trend has been initiated by the} enormous research effort in the field of distributed optimization~\cite{Boyd2011}. There exists a large number of related optimization methods, including dual decomposition~\cite{Everett1963}, ADMM~\cite{Boyd2011}, or ALADIN~\cite{Houska2016}, which have all been applied to MPC in various contexts and by many authors~\cite{Conte2012,Giselesson2013,Kozma2013,Necoara2008,Donoghue2013,Richter2011}. \change{Additionally, applications of accelerated variants of ADMM to MPC can be found in~\cite{FS16,SK13}.}

\change{However, modern distributed optimization methods, such as ADMM or ALADIN, typically converge to an optimal solution in the limit, i.e., if the number of iterations tends to infinity. Thus, if real-time constraints are present, one could at most implement a finite number of such ADMM or ALADIN iterations returning a control input that may be infeasible or sub-optimal by a large margin. But, unfortunately, for such heuristic implementations of real-time distributed MPC, there are, at the current status of research, no stability, feasibility, and performance guarantees available.}

\change{Therefore, this} paper asks the question whether it is possible to approximate an MPC feedback law by a finite code list, whose input is the current state measurement and whose output, the control reaction, is obtained by evaluating a constant, finite number of pre-computed, explicit solution maps that are associated to MPC problems of a smaller scale. Here, a key requirement is that \change{recursive feasibility, uniform asymptotic stability, and performance guarantees} of the implemented closed-loop controller have to be verifiable offline. Notice that such an MPC code would have major advantages for an embedded hardware system, \change{as it has a constant run-time using static memory only, while, at the same time, feasibility, stability, and performance guarantees are available.}

The contribution of this paper is the development of a controller, which meets these requirements under the restricting assumption that the original MPC problem is a strongly convex (but potentially large-scale) QP, as introduced in Section~\ref{sec::problem}. The control scheme itself is presented in the form of Algorithm~\ref{alg::dempc} in Section~\ref{sec::algorithm}. This algorithm alternates---similar to the distributed optimization method ALADIN~\cite{Houska2016}---between solving explicit solution maps that are associated with small-scale decoupled \change{QPs and solving a linear equation system} of a larger scale. 
However, in contrast to ALADIN, ADMM or other existing distributed optimization algorithms, Algorithm~\ref{alg::dempc} performs only a constant number of iterations \change{per sampling time.}

\change{The recursive feasibility, stability, and performance
properties} of Algorithm~\ref{alg::dempc}, which represent the main
added value compared to our preliminary work~\cite{Oravec2017}, are summarized in Section~\ref{sec::convergence},~\ref{sec::stability}, and~\ref{sec::performance}, respectively. Instead of relying on existing analysis concepts from the field of distributed optimization, the mathematical developments in this paper rely on results that find their origin in Explicit MPC theory~\cite{BorrelliPHD}. In particular, the technical developments around Theorem~\ref{thm::convergenceRate} make use of the solution properties of multi-parametric QPs in order to derive convergence rate estimates for Algorithm~\ref{alg::dempc}. Moreover, Theorem~\ref{thm::stability} establishes \change{an asymptotic stability guarantee of the presented real-time closed-loop scheme}. This result is complemented by Corollary~\ref{cor::performance}, which provides \change{bounds} on the sub-optimality of the presented control scheme.

Finally, Section~\ref{sec::pMPC} discusses implementation details with a particular emphasis on computational and storage complexity exploiting the fact that the presented scheme can be realized by using static memory only while ensuring a constant run-time. A spring-vehicle-damper benchmark is used to illustrate the performance of the proposed real-time scheme. Section~\ref{sec::conclusions} concludes the paper.

\section{Linear-Quadratic MPC}
\label{sec::problem}
This paper concerns differential-algebraic model preditive control problems in discrete-time form,\change{
\begin{eqnarray}
\begin{array}{rcl}
J(x_0) = & \underset{x,u,z}{\min} & \overset{N-1}{\underset{k=0}{\sum}} \ell(x_k,u_k,z_k) + \mathcal M( x_N ) \\[0.25cm]
& \mathrm{s.t.} & \left\{
\begin{array}{l}
\forall k \in \{ 0, \ldots, N-1 \}, \\
\begin{array}{rcl}
x_{k+1} &=& A x_k + B u_k + C z_k, \\
0 &=& D x_k  + E z_k,
\end{array} \\
x_k \in \mathbb X, \; u_k \in \mathbb U, \; z_k\in\mathbb{Z}, \; x_N \in \mathbb X_N,
\end{array}
\right.
\end{array} \label{eq::mpc}
\end{eqnarray}}
with strictly convex quadratic stage and terminal cost,
\begin{eqnarray}
\ell(x,u,z) &=& x^\tr Q x + u^\tr R u + z^\tr S z \notag \\[0.16cm]
\text{and} \qquad \change{\mathcal M(x)} &=& \change{x^\tr P x} \; . \notag
\end{eqnarray}
Here, $x_k \in \mathbb R^{n_{\mathrm{x}}}$ denotes the predicted state at time $k$, $z_k \in \mathbb R^{n_{\mathrm{z}}}$ an associated algebraic state, and $u_k \in \mathbb R^{n_{\mathrm{u}}}$ the associated control input. The matrices
\begin{align}
\begin{array}{c}
A,\change{P},Q \in \mathbb R^{n_{\mathrm{x}} \times n_{\mathrm{x}}}, \; B \in \mathbb R^{n_{\mathrm{x}} \times n_{\mathrm{u}}}, \; C \in \mathbb R^{n_{\mathrm{x}} \times n_{\mathrm{z}}},\\
\change{D \in \mathbb R^{n_{\mathrm{z}} \times n_{\mathrm{x}}}}, \; \change{E,S} \in \mathbb R^{n_{\mathrm{z}} \times n_{\mathrm{z}}}, \; R \in \mathbb R^{n_{\mathrm{u}} \times n_{\mathrm{u}}}
\end{array} \notag
\end{align}
are given and constant. Notice that~\eqref{eq::mpc} is a parametric optimization problem with respect to the current state measurement $x_0$, i.e., the optimization variable $x = [ x_1, x_2, \ldots, x_N ]$ includes all but the first element of the state sequence. The control sequence $u = [u_0, u_1, \ldots, u_{N-1}]$ and algebraic state sequence $z = [ z_0, z_1, \ldots, z_N ]$ are defined accordingly.
\begin{assumption}
\label{ass::blanket}
We assume that
\begin{enumerate}
\item[a)] the constraint sets $\mathbb U \subseteq \mathbb R^{n_{\mathrm{u}}}$, $\mathbb X, \mathbb X_N \subseteq \mathbb R^{n_{\mathrm{x}}}$, and \change{$\mathbb Z \subseteq \mathbb R^{n_{\mathrm{z}}}$} are convex and closed polyhedra satisfying $0 \in \mathbb U$, $0 \in \mathbb X$, $0 \in \mathbb X_N$, \change{and $0 \in \mathbb Z$};
        
\item[b)] the matrices $Q $, $R $, $S $, and \change{$P$} are all symmetric and positive definite;
        
\item[c)] and the \change{square-matrix $E$ is} invertible.

\end{enumerate}
\end{assumption}    
Notice that the latter assumption implies that, at least in principle, one could eliminate the algebraic states, as the algebraic constraints in~\eqref{eq::mpc} imply
\begin{equation}
\label{eq::replace}
z_k = - E^{-1} D x_k \; .
\end{equation}
However, the following algorithmic developments exploit the particular structure of~\eqref{eq::mpc}. This is relevant in the context of large-scale interconnected control \change{systems}, where the \change{matrix $D$ is} potentially dense while all other matrices have a block-diagonal structure, \change{as explained in the sections below}.

\begin{remark}
Assumption~\ref{ass::blanket}c) also implies that the equality constraints in~\eqref{eq::mpc} are linearly independent. Thus, the associated co-states (dual solutions) are unique. Similarly, Assumption~\ref{ass::blanket}a) and \ref{ass::blanket}b) imply strong convexity such that the primal solution of~\eqref{eq::mpc} is unique whenever it exists. 
\end{remark}

\StartChange
\begin{remark}
The assumption that the matrices $Q $, $R $, $S$, and $P$ are positive definite can be satisfied in practice by adding suitable positive definite regulatization to the stage cost.
\end{remark}
\EndChange

\subsection{Interconnected systems}
\label{sec::networks}
\begin{figure}[tbp!]
    \centering
    \vspace{0.5cm}
    \includegraphics[width=0.47\textwidth]{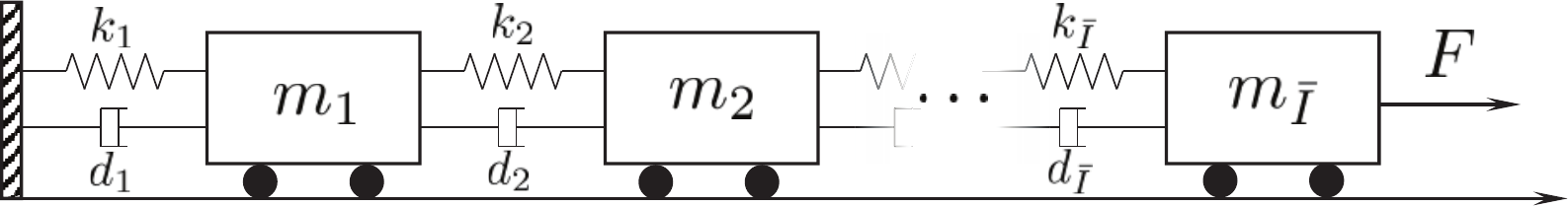}
    \caption{\change{Sketch of a spring-vehicle-damper system.}}
    \label{fig::modeling}
\end{figure}

\StartChange
Many control systems of practical interest have a particular structure, as they can be divided into $\bar I \in \mathbb N$ subsystems that are interconnected. A prototypical example for such an interconnected system is shown in Figure~\ref{fig::modeling},
which consists of $\bar I$ vehicles with mass $m$ that are linked by springs and dampers. Such systems can be modelled by a linear discrete-time recursion of the form
\begin{align}
\label{eq::network}
\underbrace{\left[ x_{k+1} \right]_i = \mathcal A_{ii} \left[ x_{k} \right]_i + \mathcal{B}_i \left[ u_k \right]_i}_{\text{dynamics of the $i$-th subsystem}}
+ \underbrace{\sum_{j \in \mathcal N_i} \mathcal A_{ij} \left[ x_{k} \right]_j}_{\substack{ \text{contribution} \\ \text{from neighbors}}
}
\end{align}
for all $k \in \{ 1, \ldots, N \}$. Here, $\left[ x_{k} \right]_i \in \mathbb X_i$ denotes the state, $\left[ u_{k} \right]_i \in \mathbb U_i$ the control input of the $i$th subsystem, and $\mathcal A_{ij}$ and $\mathcal{B}_i$ system matrices of appropriate dimension. The index $j$ runs over $\mathcal N_i$, the set of neighbors of the $i$th node in the system's connectivity graph.

\begin{example}
\label{ex::springmassdamper}
For the spring-vehicle-damper system in Figure~\ref{fig::modeling} the connectivity graph can be specified by setting
$$\mathcal N_i = \{ i-1, i+1 \}$$
for $i \in \{ 2, 3, \ldots, \bar I-1 \}$ and $\mathcal N_1 = \{ 2 \}$ and $\mathcal N_{\bar{I}}$ for the first and last trolley.
In this example, the $i$-th subblock of the state at time $k$,
$$[x_k]_i = \left( [p_k]_i \, , \, [v_k]_i \right)^\tr ,$$
consists of the position and velocity of the $i$-th trolley, all relative to their equilibrium values.
The control input is the force at the last trolley. 
A corresponding system model is then obtained by setting
\begin{align}
\begin{array}{rclrcl}
\mathcal{A}_{i,i} &=& \mathbb I + h \left(
\begin{array}{cc}
0 & 1 \\ - 2 \frac{k}{m} & -2\frac{d}{m}
\end{array}
\right) \; , \qquad & \mathcal{B}_i &=& \left( \begin{array}{c} 0 \\ 0 \end{array} \right),  \notag \\
\mathcal{A}_{\bar{I},\bar{I}} &=& \mathbb I + h \left(
\begin{array}{cc}
0 & 1 \\ - \frac{k}{m} & -\frac{d}{m}
\end{array}
\right) \; , \qquad & \mathcal{B}_{\bar{I}} &=& h \left( \begin{array}{c} 0 \\ \frac{1}{m} \end{array} \right), \notag \\
\mathcal{A}_{i-1,i} &=& \mathcal A_{i,i+1} = h \left(
\begin{array}{cc}
0 & 0  \\  \frac{k}{m} & \frac{d}{m}
\end{array}
\right),
\end{array}
\end{align}
for all $i \in \{ 1, \ldots ,\bar{I}-1 \}$. In this context, $m > 0$ denotes the mass, $k > 0$ the spring constant, $d \geq 0$ a damping coefficient, and $h>0$ the step-size of the Euler discretization.
\end{example}

System~\eqref{eq::network} can be reformulated by introducing auxiliary variables of the form \EndChange
\begin{equation}
\label{eq::refZ}
\left[ z_{k} \right]_i = \sum_{j \in \mathcal N_i} \mathcal A_{ij} \left[ x_{k} \right]_j \; .
\end{equation}
The advantage of introducing these algebraic auxiliary variables is that the dynamic recursion~\eqref{eq::network} can be written in the form
\begin{align}
x_{k+1} \;=\;& A x_k + B u_k + C z_k, \notag \\
0 \;=\;& D x_k  + E z_k, \notag
\end{align}
with $A = \mathrm{diag}\left( \mathcal A_{1,1}, \ldots, \mathcal A_{\bar{I},\bar{I}} \right)$, $B = \mathrm{diag}\left( \mathcal B_1, \ldots, \mathcal B_{\bar I} \right)$, $C = I$,
\StartChange
\[
D_{i,j} = \left\{
\begin{array}{ll}
\mathcal A_{i,j} & \text{if} \; j \in \mathcal N_i \\
0 & \text{otherwise},
\end{array}
\right\}
\]
\EndChange
and $E = -I$. After this reformulation, all matrices in the algebraic discrete-time system are block diagonal except for the matrix $D$, which may, however, still be sparse depending on the particular definition (graph structure) of the sets $\mathcal N_i$. Notice that the associated state and control constraint sets $\mathbb X = \mathbb X_1 \times \ldots \times \mathbb X_{\bar I}$ and $\mathbb U = \mathbb U_1 \times \ldots \times \mathbb U_{\bar I}$ have a separable structure, too.

\subsection{Recursive feasibility and asymptotic stability}

Notice that the stability and recursive feasibility properties of MPC \change{controllers have} been analyzed exhaustively~\cite{Rawlings2009}. \change{
As long as Assumption~\ref{ass::blanket} holds, these results can be applied one-to-one for~\eqref{eq::mpc} after eliminating the algebraic states explicitly.}
\StartChange
\begin{definition}
A set $X \subseteq \mathbb X$ is called control invariant if
\[
X \subseteq \left\{ x \in \mathbb R^{n_{\mathrm{x}}} \left|
\begin{array}{l}
\exists x^+ \in X, \, \exists u \in \mathbb U, \, \exists z \in \mathbb Z, \\
\begin{array}{rcl}
x^+ &=& A x + B u + C z, \\
0 &=& D x + E z
\end{array}
\end{array}
\right.
\right\} \; .
\]
\end{definition}
\EndChange
\begin{assumption}
\label{ass::forwardInvariance}
We assume that the terminal \change{set $\mathbb X_N$ is control invariant.}
\end{assumption}
It is well known that~\eqref{eq::mpc} is recursively feasible if Assumption~\ref{ass::blanket} and Assumption~\ref{ass::forwardInvariance} hold \change{recalling that one can always eliminate the algebraic states in order to apply the results from~\cite{Rawlings2009}}.
Another standard assumption can be formulated as follows.
\begin{assumption}
\label{ass::TerminalLyapunov}
The terminal cost $\mathcal M$ in~\eqref{eq::mpc} admits a control law $\mu: \mathbb X_N \to \mathbb U$ such that for all $x \in \mathbb X_N$
\[
\change{\ell(x,\mu(x),z) + \mathcal M( x^+ ) \leq \mathcal M(x)} \; ,
\]
where $z = -E^{-1} D x$ \change{and $x^+ = A x + B \mu(x) - C E^{-1} D x$}. 
\end{assumption}
The MPC controller~\eqref{eq::mpc} is asymptotically stable if Assumption~\ref{ass::blanket},~\ref{ass::forwardInvariance}, and~\ref{ass::TerminalLyapunov} hold~\cite{Rawlings2009}. Notice that there exist \change{generic} methods for both the construction of forward invariant sets $\mathbb X_N$ and the construction of quadratic terminal costs such that the above criteria are satisfied~\cite{Bla:aut:99,Rawlings2009}. \change{Moreover, tailored methods for constructing block separable terminal costs and terminal constraint sets are available in the literature, too~\cite{Conte2016}.}

\section{Suboptimal real-time MPC}
\label{sec::algorithm}
In this section we propose and analyze a real-time algorithm for finding approximate, suboptimal solutions of~\eqref{eq::mpc}.

\subsection{Preliminaries}
Let us introduce the stacked \change{vectors
\[
y_0 = \left[ u_0^\tr \, , \, z_0^\tr \right]^\tr
\; , \; \;
y_k = \left[ x_k^\tr \, , \, u_k^\tr \, , \, z_k^\tr
\right]^\tr \; , \; \; 
y_N = x_N,
\]
and their associated constraint sets
\begin{eqnarray}
\mathbb Y_k &= \left\{ y_k
\left|
\begin{array}{l}
A x_k + B u_k + C z_k \in \mathbb X, \\
0 = D x_k + E z_k, \\
x_k \in \mathbb X, \; u_k \in \mathbb U, \; z_k \in \mathbb Z
\end{array}
\right.
\right\}, \;
\mathbb Y_N &= \mathbb X_N,
\end{eqnarray}
for all $k \in \{ 0, \ldots, N-1 \}$.} Moreover, it is convenient to introduce the shorthand notation
\begin{equation}
\label{eq:def_F}
\begin{split}
F_k(y_k) &\;=\; \ell(x_k,u_k,z_k) \; , \quad k\in \{ 0, \ldots, N-1 \} \\[0.1cm]
F_N(y_N) &\;=\; \mathcal M(x_N) , 
\end{split}
\end{equation}
as well as the matrices
\[
G_k = \left(
\begin{array}{ccc}
I & 0 & 0 \\
0 & 0 & 0
\end{array}
\right) \; , \; \change{G_{N} = I},
\]
for $k \in \{ 1, \ldots N-1 \}$ and
\[
H_0 = \left(
\begin{array}{ccc}
B & C \\
0 & E
\end{array}
\right)
, \,
H_k = \left(
\begin{array}{ccc}
A & B & C \\
D & 0 & E
\end{array}
\right),
\]
as well as $h_0 = [ (A x_0)^\tr \, 0]^\tr, \;  h_{k} = 0$ for all indices $k \in \{ 1, \ldots, N-1 \}$. The advantage of this notation is that~\eqref{eq::mpc} can be written in the \change{form
\begin{align}
\label{eq::mpc2}
J(x_0) = \underset{y}{\min} &\quad \sum_{k=0}^{N} F_k( y_k ) \\[0.16cm] \notag
\text{s.t.} &\quad 
\left\{
\begin{array}{l}
\forall k \in \{ 0, \ldots, N-1 \}, \\
\begin{array}{rcll}
G_{k+1} y_{k+1} &=& H_k y_k + h_k \; & \mid \; \lambda_{k},
\end{array} \\
y_k \in \mathbb Y_k \; , \; y_N \in \mathbb Y_N \; .
\end{array}
\right.
\end{align}
Here, $\lambda_0, \ldots, \lambda_{N-1}$} denote the multipliers of the affine constraints in~\eqref{eq::mpc2}. It is helful to keep in mind that both the function $F_0$ and the vector $h_0$ depend on $x_0$ even if we hide this parametric dependence for a moment in order to arrive at a simpler notation. In addition, we introduce a shorthand for the objective in~\eqref{eq::mpc2},
\begin{align}
F(y) =& \sum_{k=0}^{N} F_k( y_k ) \; . \notag
\end{align}
Moreover, the convex conjugate function of $F$ is denoted by
\begin{align}
F\opt(\lambda) =& \max_{y}\;\ \bigg\{ - F(y) - \left( H_0^\tr \lambda_{0}^m \right)^\tr y_0  \notag \\
& + \sum_{k=1}^{N-1} \left( G_k^\tr \lambda_{k-1}^m - H_k^\tr \lambda_{k}^m \right)^\tr y_k +  \lambda_{N-1}^\tr G_N y_N
\bigg\} \; . \notag
\end{align}
Notice that the functions $F$ and $F\opt$ are strongly convex quadratic forms with $F(0) = 0$ and $F\opt(0)=0$ as long as Assumption~\ref{ass::blanket} is satisfied. The optimal primal and dual solutions of~\eqref{eq::mpc2} are denoted by $x\opt$ and $\lambda\opt$, respectively. It is well-known that $x\opt$ and $\lambda\opt$ are continuous and piecewise affine functions of $x_0$ on the polytopic domain $\mathcal X = \{ x_0 \mid  J(x_0) < \infty \}$, see~\cite{borrelli2003geometric}.

\subsection{Algorithm}
\StartChange
The main idea for solving~\eqref{eq::mpc2} approximately and in real time is to consider the unconstrained auxiliary optimization problem
\begin{align}
\label{eq::mpc2mod}
J(x_0) = \underset{y}{\min} &\quad \sum_{k=1}^{N-1} F_k( y_k - y_k^\text{ref} ) + F_N(y_N- y_N^\text{ref})\\[0.16cm] \notag
\text{s.t.} &\quad 
\left\{
\begin{array}{l}
\forall k \in \{ 0, \ldots, N-1 \}, \\
\begin{array}{rcll}
G_{k+1} y_{k+1} &=& H_k y_k + h_k \; & \mid \; \lambda_{k},
\end{array}
\end{array}
\right.
\end{align}
which is a standard tracking optimal control problem without inequality constraints. Here, we have introduced the reference trajectory $y_k^\text{ref}$. For the special case that $y^\text{ref} = y\opt$ is equal to the minimizer of~\eqref{eq::mpc2}, Problems~\eqref{eq::mpc2} and~\eqref{eq::mpc2mod} are equivalent. Notice that the main motivation for introducing the coupled QP~\eqref{eq::mpc2mod} is that this problem approximates~\eqref{eq::mpc2} without needing inequality constraints. Thus, this problem can be solved by using a sparse linear algebra solver; see Section~\ref{sec::sparseSolver} for implementation details.

\begin{algorithm}[htbp!]
    \caption{Parallel real-time MPC}
    \textbf{Initialization:}\\[0.2cm]
    Choose initial $y^1 = [y_0^1, \ldots, y_N^1]$ and \change{$\lambda^1 = [ \lambda_0^1, \ldots, \lambda_{N-1}^1 ]$} and a constant $\gamma > 0$.\\[0.2cm]
    \textbf{Online:}\\
    \begin{enumerate}
        \item[\textbf{1)}] Wait for the state measurement $x_0$ and compute the constant
        \[
        f^1 = F(y^1) + F\opt(\lambda^1) \; .
        \]
        If $f^1 \geq \gamma^2 x_0^\tr Q x_0$, rescale
        $$y^1 \leftarrow y^1 \sqrt{\frac{\gamma^2 \Vert x_0 \Vert_Q^2}{f_1}} \; \; \text{and} \; \; \lambda^1 \leftarrow \lambda^1 \sqrt{\frac{\gamma^2 \Vert x_0 \Vert_Q^2}{f_1}} \; ,$$ \\[0.1cm]
        where $\Vert x_0 \Vert_Q^2 \triangleq x_0^{\tr} Q x_0 $ .
        \item[\textbf{2)}] \textbf{For $m = 1\rightarrow \overline{m}$}\\[0.2cm]
        \noindent
        \textbf{a)} solve the small-scale decoupled QPs in parallel
        \begin{eqnarray}
        \hspace{-0.6cm} &\underset{\xi_0^m \in \mathbb Y_0}{\min} & \hspace{-0.1cm}  F_0(\xi_0^m)  -  (H_0^\tr \lambda_{0}^m)^\tr \xi_0^m + F_0( \xi_0^m - y_0^m )  \notag \\[0.2cm]
        \hspace{-0.6cm} &\underset{\xi_k^m \in \mathbb Y_k}{\min} & \hspace{-0.1cm} F_k(\xi_k^m)  + ( G_k^\tr \lambda_{k-1}^m - H_k^\tr \lambda_{k}^m )^\tr \xi_k^m + F_k( \xi_k^m - y_k^m )   \notag \\[0.2cm]
        \hspace{-0.6cm} &\change{ \underset{\xi_N^m \in \mathbb Y_N}{\min}} & \change{\hspace{-0.1cm} F_N(\xi_N^m)  + \left( G_N^\tr \lambda_{N-1}^m \right)^\tr \xi_N^m + F_N( \xi_N^m - y_N^m ) }  \notag
        \end{eqnarray}
        for all \change{$k \in \{ 1, \ldots, N-1 \}$} and denote the optimal solutions by $\xi^m = [ \xi_0^m, \xi_1^m, \ldots, \xi_N^m ]$.\\[0.2cm]
        \textbf{b)} Solve the coupled QP
        \begin{align}
        \label{eq::qp}
        \underset{y}{\min}& \quad \sum_{k=0}^{N} F_k( y_k^{m+1} - 2 \xi_k^m + y_k^m ) \\[0.16cm]\notag
        \text{s.t.} &\quad 
        \left\{
        \begin{array}{l}
        \forall k \in \{ 0, \ldots, N-1 \}, \\
        \change{
            \begin{array}{rcll}
            G_{k+1} y_{k+1}^{m+1} &=& H_k y_k^{m+1} + h_k \; & \mid \; \delta_{k}^m,
            \end{array}
        }
        \end{array}
        \right.
        \end{align}
        and set $\lambda^{m+1} = \lambda^m + \delta^m$.\\[0.2cm]
        \noindent
        \textbf{End}\\
        \item[\textbf{3)}] Send $u_0 = \left( 0 \;\; I \;\; 0  \right) \xi_0^{\overline{m}}$ to the real process.\\
        
        \item[\textbf{4)}] Set $y^1 = [y_1^{\overline{m}}, \ldots, y_N^{\overline{m}}, 0]$\change{, $\lambda^1 = [ \lambda_1^{\overline{m}}, \ldots, \lambda_{N-1}^{\overline{m}}, 0]$} and go to Step~1.\\
    \end{enumerate}
    \label{alg::dempc}
\end{algorithm}

Let us assume that $y^m$ and $\lambda^m$ are the current approximations of the primal and dual solution of~\eqref{eq::mpc2}. Now, the main idea of Algorithm~\ref{alg::dempc} to construct the next iterate $y^{m+1}$ and $\lambda^{m+1}$ by performing two main operations. First, we
solve augmented Lagrangian optimization problems of the form
\begin{align}
\label{eq::NLPs}
\underset{\xi^m \in \mathbb Y}{\min} & \;\; F(\xi^m)  -  (\mathcal G^\tr \lambda_{0}^m)^\tr \xi^m + \left\| \xi^m - y^m \right\|_{\mathcal Q}^2
\end{align}
with
$$\mathcal G = \left(
\begin{array}{ccccc}
H_0 & -G_0 &  & & 0 \\
 & H_1 & - G_1 & & \\
 & & \ddots & \ddots & \\
0 & & & H_{N-1} & -G_{N-1}
\end{array}
\right).
$$
The focus of the following analysis is on the case that the weighting matrix $\mathcal Q = \nabla^2 F(0)$ is such that
$$\left\| \xi^m - y^m \right\|_{\mathcal Q}^2 = F( \xi^m - y^m )$$
recalling that $F$ is a centered positive-definite quadratic form. And second, we solve QP~\eqref{eq::mpc2mod} for the reference
point
\[
y^\text{ref} = 2 \xi^m - y^m \; ,
\]
which can be interpreted as weighted average of the solution of~\eqref{eq::NLPs} and the previous iterate $y^m$. These two main steps correspond to Step 2a) and Step 2b) in  Algorithm~\ref{alg::dempc}. Notice that the main motivation for introducing the augmented Lagrangian problem~\eqref{eq::NLPs} is that this optimization problem is separable, as exploited by Step~2a). As explained in more detail in Section~\ref{sec::parametric}, the associated smaller-scale QPs can be solved by using existing tools from the field of Explicit MPC.

Additionally, in order to arrive at a practical procedure, Algorithm~\ref{alg::dempc} is initialized with guesses,
$$y^1 = [y_0^1, \ldots, y_N^1] \quad \text{and} \quad \lambda^1 = [ \lambda_0^1, \ldots, \lambda_N^1 ] \; ,$$
for the primal and dual solution of~\eqref{eq::mpc2}. Notice that Algorithm~\ref{alg::dempc} receives a state measurement $x_0$ in every iteration (Step 1) and returns a control input to the real process (Step 3).
Moreover, Step 1) rescales $y^1$ and $\lambda^1$ based on a tuning parameter $\gamma > 0$, which is assumed to have the following property.
\begin{assumption}
\label{ass::gamma}
The constant $\gamma$ in Algorithm~\ref{alg::dempc} is such that
\[
F(y\opt) + F\opt(\lambda\opt) \leq \gamma^2 x_0^\tr Q x_0 \; .
\]
\end{assumption}
Notice that such a constant $\gamma$ exists and can be pre-computed offline, which simply follows from the fact that $y\opt$ and $\lambda\opt$ are Lipschitz continuous and piecewise affine functions of $x_0$, as established in~\cite{borrelli2003geometric}.
The rationale behind this rescaling is that this step prevents initializations that are too far away from $0$. Intuitively, if the term $f^1 = F(y^1) + F\opt(\lambda^1)$ is much larger than $x_0^\tr Q x_0$, then $(y^1,\lambda^1)$ can be expected to be far away from the optimal solution $(y^\star,\lambda^\star)$ and it is better to rescale these variables such that they have a reasonable order of magnitude. In the following section, we will provide a theoretical justification of the rescaling factor $\sqrt{\frac{\gamma^2 \Vert x_0 \Vert_Q^2}{f_1}}$. Notice that if we would set this rescaling factor to $1$, Algorithm 1 is unstable in general. In order to see this, consider the scenario that a user initializes the algorithm with an arbitrary $(y^1,\lambda^1) \neq 0$. Now, if the first measurement happens to be at $x_0 =0$, of course, the optimal control input is at $u^\star =0$. But, if we run Algorithm~1 with $\bar m < \infty$, it returns an approximation $u_0 \approx u^\star = 0$, which will introduce an excitation as we have $u_0 \neq 0$ in general. Thus, if we would not rescale the initialization in Step~1), it would not be possible to establish stability.
\EndChange

\subsection{Convergence properties of Algorithm~\ref{alg::dempc}}
\label{sec::convergence}

This section provides a concise overview of the theoretical convergence properties of Algorithm~\ref{alg::dempc}. We start with the following theorem, which is proven in Appendix~\ref{app::convergenceRate}.
\begin{theorem}
\label{thm::convergenceRate}
\change{Let Assumption~\ref{ass::blanket} be satisfied and let~\eqref{eq::mpc2} be feasible, i.e., such that a unique minimizer $y\opt$ and an associated dual solution $\lambda\opt$ exist.}  Then there exists a positive constant $\kappa < 1$ such that
\begin{align}
&\;\; F( y^{m+1} - y\opt ) + F\opt( \lambda^{m+1} - \lambda\opt ) \notag \\[0.16cm]
\leq& \;\;\kappa \left( F( y^m - y\opt ) + F\opt( \lambda^m - \lambda\opt ) \right)
\end{align}
for all $m \geq 2$.
\end{theorem}
\change{Because Theorem~\ref{thm::convergenceRate} establishes contraction, an immediate consequence is that the iterates of Algorithm~1 would converge to the exact solution of~\eqref{eq::mpc2}, if we would set $\overline{m} = \infty$, i.e.,
\[
\lim_{m \to \infty} \xi^m = y\opt \quad \text{and} \quad  \lim_{m \to \infty} \lambda^m = \lambda\opt
\]
Moreover, the proof of Theorem~\ref{thm::convergenceRate} provides} an explicit procedure for computing the constant $\kappa < 1$.


\subsection{Recursive feasibility and asymptotic stability of Algorithm~\ref{alg::dempc}}
\label{sec::stability}

The goal of this section is to establish recursive feasibility and asymptotic stability of Algorithm~\ref{alg::dempc} on $\mathcal X$. Because we send the control input $u_0 = \left(  0 \; I \; 0  \right) \xi_0^{\overline{m}}$ to the real process, the next measurement will be at
\[
x_0^+ = A x_0 + (B \;\; C) \xi_0^{\overline{m}} \; .
\]
Notice that, in general, we may have
\[
x_0^+ \neq x_1\opt = A x_0 + (B \;\; C) y_0\opt \; ,\
\]
since we run Algorithm~\ref{alg::dempc} with a finite $\overline{m} < \infty$. \change{A proof of the following proposition can be found in Appendix~\ref{app::feasibility}.}
\StartChange
\begin{proposition}
\label{prop::feasibility}
Let us assume that
\begin{enumerate}
\item Assumptions~\ref{ass::blanket} is satisfied,
\item we use the terminal region $\mathbb X_N = \mathbb X$, and
\item $\mathbb X$ is control invariant.
\end{enumerate}
Then Algorithm~\ref{alg::dempc} is recursively feasible in the sense that $x_0 \in \mathbb X$
implies $x_0^+ \in \mathbb X$. Moreover, the equation $\mathcal X = \mathbb X$ holds, i.e.,
Problem~\eqref{eq::mpc2} remains feasible.
\end{proposition}

\begin{remark}
The construction of control invariant sets for interconnected systems in the presence
of state constraints can be a challenging task
and this manuscript does not claim to resolve this problem. An in-depth discussion of how
to meet the third requirement of Proposition~\ref{prop::feasibility} is beyond the scope of this
paper, but we refer to~\cite{Conte2016}, where a discussion of methods for the construction
of control invariant sets for network systems can be found.
\end{remark}

The following theorem establishes asymptotical stability of Algorithm~\ref{alg::dempc},
one of the main contributions of this paper.
The corresponding proof can be found in Appendix~\ref{app::stability}.
\begin{theorem}
\label{thm::stability}
Let Assumptions~\ref{ass::blanket},~\ref{ass::forwardInvariance},~\ref{ass::TerminalLyapunov}, and~\ref{ass::gamma} be satisfied and let $\mathbb X_N = \mathbb X$ be control invariant.
Let the constant $\sigma > 0$ be such that the semi-definite inequality
\[
\left(
\begin{array}{cc}
B^\tr Q B  & B^\tr Q C \\
C^\tr Q B  & C^\tr Q B
\end{array}
\right)
\preceq \sigma \left(
\begin{array}{cc}
R & 0 \\0 & S
\end{array}
\right)
\]
holds and let the constants $\eta,\tau > 0$ be such that
\begin{align}
\label{eq::etaTauBound}
|J(x_0^+) - J(x_1\opt)| \leq \eta \| x_0^+ - x_1\opt \|_Q + \frac{\tau}{2} \| x_0^+ - x_1\opt \|_Q^2
\end{align}
If the constant $\overline{m} \in \mathbb N$ satisfies
\begin{align}
\label{eq::mbarBound}
\overline{m} > 2 \frac{\log \left( \eta \sqrt{\sigma} \gamma (1+\sqrt{\kappa}) + \frac{\tau \sigma \gamma^2 (1+\sqrt{\kappa})^2}{2}   \right)}{\log(1/\kappa)} \; ,
\end{align}
then the controller in Algorithm~\ref{alg::dempc} is \change{asymptotically stable} on $\mathcal X$.
\end{theorem}
\EndChange

\change{The} constants $\eta, \tau, \sigma, \change{\gamma}$, and $\kappa < 1$ in the above theorem depend on the problem data only\change{,
but they are independent of the initial state $x_0 \in \mathcal X$.}

\StartChange
\begin{remark}
The lower bound~\eqref{eq::mbarBound} is monotonically increasing in $\eta, \tau, \sigma, \gamma$, and $\kappa$. Thus, the smaller these constants are, the tighter the lower bound~\eqref{eq::mbarBound} will be. For small-scale applications, one can compute these constants offline, by using methods from the field of explicit MPC~\cite{BemEtal:aut:02,borrelli2003geometric}. However, for large-scale applications, the explicit solution map cannot be computed with reasonable effort, not even offline. In this case, one has to fall back to using conservative bounds. For example, the constants $\eta$ and
$\tau$ satisfying~\eqref{eq::etaTauBound} can be found by using methods from the field of approximate dynamic programming~\cite{Wang2015}. However, if one uses such conservative bounding methods, the lower bound~\eqref{eq::mbarBound} is conservative, too.
\end{remark}
\EndChange

\subsection{Performance of Algorithm~\ref{alg::dempc}}
\label{sec::performance}
The result of Theorem~\ref{thm::stability} can be extended in order to derive an a-priori verifiable upper bound on the sub-optimality of Algorithm~1.

\begin{corollary}
\label{cor::performance}
Let the assumption of Theorem~\ref{thm::stability} hold with
\[
\alpha = 1 - \left[ \eta \sqrt{\sigma} \gamma (1+\sqrt{\kappa}) + \frac{\tau \sigma \gamma^2 (1+\sqrt{\kappa})^2}{2} \right] \kappa^{\frac{\overline{m}}{2}} \; .
\]
If $y_i^{\mathrm{cl}} = \left( x_i^\text{cl}, \, u_i^\text{cl}, \, z_i^\text{cl} \right)$ denotes the sequence of closed-loop states and controls that are generated by the controller in Algorithm~\ref{alg::dempc}, an a-priori bound on the associated infinite-horizon closed-loop performance is given by
\[
\sum_{i=0}^{\infty} \ell( x_i^\text{cl}, \, u_i^\text{cl}, \, z_i^\text{cl} ) \leq \frac{J(x_0)}{\alpha} \; .
\]
\end{corollary}

\noindent
\textbf{Proof.}
Theorem~\ref{thm::stability} together with~\eqref{eq::lyapD1} and~\eqref{eq::lyapD2} imply that
\[
J( x_{i+1}^\text{cl} ) \leq J( x_i^\text{cl} ) - \alpha F_0( y_i^{\mathrm{cl}} ) \; ,
\]
which yields the inequality
\[
\sum_{i=0}^{\infty} F_0( y_i^{\mathrm{cl}} ) \leq  \frac{1}{\alpha} \sum_{i=0}^{\infty} \left( J(x_{i}^\text{cl}) - J(x_{i+1}^\text{cl}) \right) \; .
\]
The statement of the corollary follows after simplifying the telescoping sum on the right and substituting the equation $F_0( y_i^{\mathrm{cl}} ) = \ell( x_i^\text{cl}, \, u_i^\text{cl}, \, z_i^\text{cl} )$.
\hfill\hfill\qed


\section{Implementation on hardware with limited memory}
\label{sec::pMPC}
This section discusses implementation details for Algorithm~\ref{alg::dempc} with a particular emphasis on run-time aspects and limited memory requirements, as needed for the implementation of MPC on embedded hardware. Here, the implementation of Steps~1),~3), and~4) turns out to be straightforward, as both the CPU time requirements and memory needed for implementing these steps is negligible compared to Step 2). Thus, the focus of the following subsections is on the implementation of Step~2a) and Step~2b).

\subsection{Parametric quadratic programming}
\label{sec::parametric}
In Step~2a) of Algorithm~\ref{alg::dempc} decoupled QPs have to be
solved on-line. To diminish the induced implementation effort, we
propose to solve these QPs off-line using multi-parametric programming, i.e., to pre-compute the solution maps
\begin{equation}
\label{eq::mQP}
\begin{array}{rcl}
\xi_0 \opt(\theta_0,x_0) &=& \underset{\xi_0 \in \mathbb{Y}_0}{\text{argmin}} \; 2 F_0(\xi_0) + \theta_0^\tr \xi_0 , \\
\xi_1 \opt(\theta_1) &=& \underset{\xi_1 \in \mathbb{Y}_0}{\text{argmin}}  \; 2 F_1(\xi_1) + \theta_1^\tr \xi_1 , \\
\xi_N \opt(\theta_N) &=& \underset{\xi_N \in \mathbb{Y}_N}{\text{argmin}} \; 2F_N(\xi_N) + \theta_N^\tr \xi_N ,
\end{array}
\end{equation}
with parameters $\theta_0 \in \mathbb R^{n_{\mathrm{u}} + n_{\mathrm{z}}}$, $\theta_1 \in \mathbb R^{n_{\mathrm{x}}+ n_{\mathrm{u}} + n_{\mathrm{z}}}$, and \change{$\theta_N \in \mathbb R^{n_{\mathrm{x}}}$}.
Because these QPs are strongly convex, the functions \change{$\xi_0 \opt$, $\xi_1 \opt$, and $\xi_N \opt$} are piecewise affine~\cite{BorrelliPHD}. Here, it should be noted that $\xi_0 \opt$ additionally depends on $x_0$ recalling that this dependency had been hidden in our definition of $F_0$ and $\mathbb{Y}_0$. In this paper, we use MPT~\cite{MPT3} to pre-compute and store the explicit solution maps $\xi_0 \opt$, $\xi_1 \opt$ and $\xi_N \opt$. Consequently, Step~2a) in Algorithm~1 can be replaced by:

\begin{quote}
\textbf{Step 2a')} Compute the parameters
\begin{subequations}
\begin{align}
\theta^m_0 &= -H_0^\tr \lambda_{0}^m - 2\Sigma_0y_0^m , \\[0.1cm]
\theta^m_k &=  G_k^\tr \lambda_{k-1}^m-H_k^\tr \lambda_{k}^m  - 2\Sigma_k y_k^m , \\[0.1cm]
\change{\theta^m_N} &= \change{ G_N^\tr \lambda_{N-1}^m  - 2\Sigma_N y_N^m }
\end{align}
\end{subequations}
with $\Sigma_0 = \text{blkdiag}\{R,S\}$, $\Sigma_k = \text{blkdiag}\{Q,R,S\}$ for all $k \in \{ 1, \ldots, N-1\}$, $\Sigma_N = \change{P}$ and set
\begin{align}
\xi_0^m = \xi_0\opt( \theta_0^m, x_0 )\;,\;\;\xi_k^m = \xi_1\opt( \theta_k^m ) \notag
\end{align}
for all $k \in \{ 1, \ldots, N \}$ by evaluating the respective
explicit solution maps~\eqref{eq::mQP}. \change{In this paper, we use the enumeration-based multi-parametric QP algorithm from~\cite{HJ15} for generating these maps.}
\end{quote}

Notice that the complexity of pre-processing the small-scale QPs~\eqref{eq::mQP} solely depends on the maximum number $N_{\text{R}} = \max \{  N_{\text{R},0}, N_{\text{R},1}, N_{\text{R},N} \}$ of critical regions over which the PWA optimizers $\xi_0 \opt$, $\xi_1 \opt$ and $\xi_N \opt$ are defined~\cite{BemEtal:aut:02}, but $N_{\mathrm{R}}$ is independent of the prediction horizon $N$ as summarized in Table~\ref{tab::complexity}; see also~\cite{Oravec2017,BoyVan:ConOpt:04,Borrelli2017}. 
Here, we assume that each parametric QP is post-processed, off-line,
to obtain binary search trees~\cite{TJB03} in $\mathcal{O}(N_R^2)$ time. Once the
trees are constructed, they provide for a fast evaluation of the solution
maps in~\eqref{eq::mQP} in time that is logarithmic in the number of
regions, thus establishing the $\mathcal{O}(N \log_2(N_R))$ on-line
computational bound. The memory requirements are directly proportional
to the number of regions $N_R$ with each region represented by a
finite number of affine half-spaces.

%

\begin{table}[htbp!]
\caption{\label{tab::complexity} Computational and storage complexity of Steps~2a') and 2b) of Algorithm~\ref{alg::dempc}.}
\begin{center}
\begin{tabular}{|l|c|c|c|}
\hline
Step & Offline & Online & Memory \\
 & CPU time & CPU time & Requirement \\
\hline
2a') & $\mathcal{O}(N_{\mathrm{R}}^2)$ & $\mathcal{O}(N \log_2( N_{\mathrm{R}}))$ & $\mathcal{O}(N_{\mathrm{R}})$ \\
2b) & $\mathcal{O}( N n^3 )$ & $\mathcal{O}( N n^2 )$ & $\mathcal{O}( N n^2 )$\\
\hline
\end{tabular}
\end{center}
\bigskip
\end{table}

\subsection{Sparse linear algebra}
\label{sec::sparseSolver}
In Step~2b) of Algorithm~\ref{alg::dempc} the large-scale, coupled QP~\eqref{eq::qp} must be solved. Because this QP has equality constraints only,~\eqref{eq::qp} is equivalent to a large but sparse \change{system of equations}. Moreover, all matrices in~\eqref{eq::qp} are given and constant during the online iterations. This means that all linear algebra decompositions can be pre-computed offline. If one uses standard Riccati recursions for exploiting the band-structure of~\eqref{eq::qp}, the computational complexity for all offline computations is at most of order $\mathcal{O}( N n^3 )$, where $n = n_{\mathrm{x}} + n_{\mathrm{z}}$, while the online implementation has complexity $\mathcal{O}( N n^2 )$~\cite{Bertsekas2012}. If one considers interconnected systems this run-time result may be improved---in many practical cases, e.g., in the example that is discussed in Section~\ref{sec::caseStudy}, one may set \change{$n = ( n_{\mathrm{x}} + n_{\mathrm{z}} ) / \bar I$}. However, in general, the choice of the linear algebra solver and its computational complexity depend on the particular topology of the network~\cite{Borrelli2017}.

In summary, Algorithm~\ref{alg::dempc} can be implemented by using static memory only allocating at most $\mathcal{O}(N_{\mathrm{R}} + N n^2)$ floating point numbers. Here, the explicit solutions maps of both the decoupled QPs in Step~2a) as well as the coupled QP in Step~2b) can be pre-computed offline. Because Theorem~\ref{thm::stability} provides an explicit formula for computing a constant number of iterations $\overline{m}$ such that a stable and recursively feasible controller is obtained, the online run-time of Algorithm~\ref{alg::dempc} is constant and of order $\mathcal{O}(N \log_2( N_{\mathrm{R}}) + N n^2)$. Thus, Algorithm~\ref{alg::dempc} may be called an explicit MPC method---in the sense that it has a constant run-time and constant memory requirements for any given $N$ and $\bar I$ while stability and recursive feasibility can be verified offline. Its main advantage compared to existing MPC controllers is that it scales up easily for large scale interconnected networks of systems as well as long prediction horizons $N$, as illustrated below.

\section{Numerical example}
\label{sec::caseStudy}

\change{This section applies Algorithm~\ref{alg::dempc} to a spring-vehicle-damper control system,
which has been introduced in Example~\ref{ex::springmassdamper} with
state and control and constraints}
\begin{align*}
& \mathbb X = \mathbb X_1 \times \ldots \times \mathbb X_{\bar I} \; , \; \mathbb U = \change{[-2,0.5]} \;,\; \mathbb Z = \mathbb R^{2 \bar I} \; , \\[0.16cm]
& \text{where} \quad \mathbb X_1 = \ldots = \mathbb X_{\bar I} = \change{[-0.5,1.5]} \times \change{[-0.5,1]} \; .
\end{align*}
\change{The} weighting matrices of the MPC objective are set to
$$Q= 10 \,  I, \; R = I, \quad \text{and} \quad S = 10^{-2} \, I \; .$$
The numerical values for the mass $m$, spring constant $k$, and damping constant $d$ are listed below.
\begin{center}
\scriptsize
\begin{tabular}{lcc}
    \toprule
    Parameter & Symbol & Value\,[Unit] \\
    \midrule
    sampling time & $h$ &$0.1$\,[s]\\[0.1cm]
    spring constant & $k$ & $3$\,[N/m] \\[0.1cm]
    mass of vehicle & $m$ & $1$\,[kg]  \\[0.1cm]
    viscous damping coefficient & $d$ & $3$\,[N\,s/m]\\[0.1cm]
    \bottomrule
\end{tabular}
\end{center}
Last but not least the \change{matrix $P$ is computed} by solving an algebraic Riccati equation, such that the terminal cost is equal to the unconstrained infinite horizon cost if none of the inequality constraints are \change{active~\cite{Rawlings2009}}.

\change{We have implemented Algorithm~\ref{alg::dempc} in Matlab R2018a using YALMIP~\cite{yalmip} and MPT 3.1.5~\cite{MPT3}. Here, the solution maps of the QPs~\eqref{eq::mQP} were pre-computed using the geometric parametric LCP solver of MPT 3.1.5~\cite{MPT3}}. By exploiting the separability of the cost function and constraints, the storage space for the parametric solutions can be reduced to $287 \, \text{kB}$, which corresponds to $432$ critical regions. This memory requirement is independent of the number of vehicles $\bar{I}$ and the length of the prediction horizon $N$. In contrast to this, the number of regions for standard Explicit MPC depends on both $\bar{I}$ and $N$:
\begin{center}
\scriptsize
\begin{tabular}{c | c c}
    $(\bar{I},N)$ & \# of regions  & memory\\
    \hline
    $(1,10)$ &  $ 58$ & $14 \,  [\text{kB}]$ \\
    $(1,20)$ &  $84$ & $ 40 \,  [\text{kB}] $ \\
    $(1,50)$ &  $144$ & $ 169 \,  [\text{kB}] $ \\
    $(2,10)$ &  $2244$ &  $ 877 \,  [\text{kB}] $ \\
    $(3,10)$ & $4247$  &$2324 \,  [\text{kB}]$ \\
\end{tabular}
\end{center}
Notice that the number of regions of \change{standard~Explicit MPC} explode quickly, as soon as one attempts to choose $N \geq 10$ or more than $3$ vehicles. In contrast to this, our approach scales up easily to hundreds of vehicles and very long horizons. Here, only the evaluation of the solution map of~\eqref{eq::qp} depends on $\bar{I}$ and $N$. For example, if we set $\bar{I} = 3$ and $N=30$, we need $10 \, \text{kB}$ to store this map---for larger values this number scales up precisely as predicted by Table~\ref{tab::complexity}.

\begin{figure}[htbp!]
    \centering
    \includegraphics[width=0.4\textwidth]{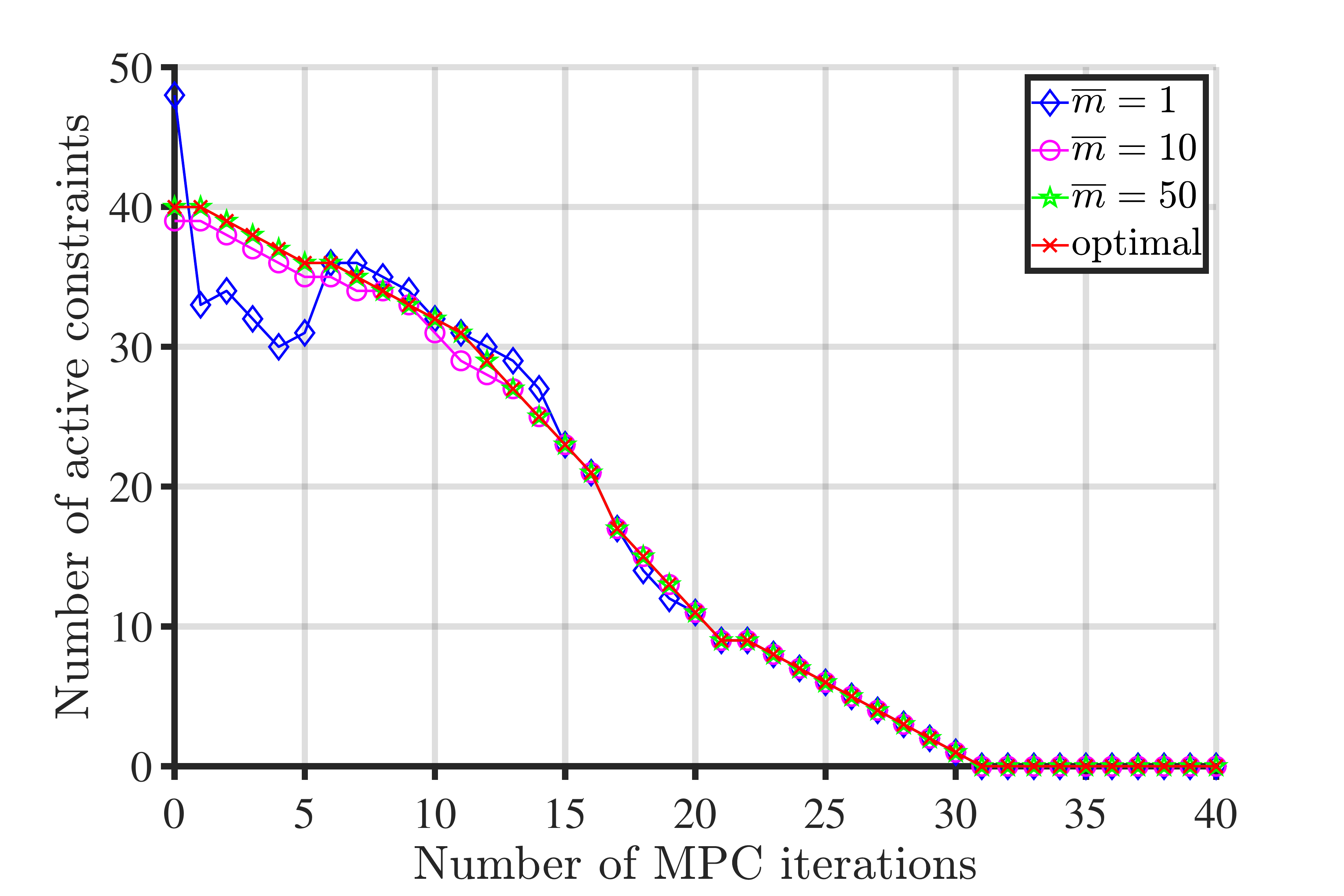}
    \caption{\change{The total number of active constraints of all distributed QP solvers during the MPC iterations for different choices of $\overline{m}$.}}
    \label{fig::active_set}
\end{figure}

\change{Figure~\ref{fig::active_set} shows the total number of active constraints of all distributed QP solvers during the MPC iterations for different choices of $\overline{m}$. In order to visualize the performance of the proposed sub-optimal controller in terms of the number of correctly detected active constraint indices, the number of active constraints of non-suboptimal MPC (corresponding to $\overline{m} = \infty$) are shown in the form of red crosses in Figure~\ref{fig::active_set}. If we compare these optimal red crosses with the blue diamonds, which are obtained for $\overline{m} = 1$, we can see that the choice $\overline{m} = 1$ still leads to many wrongly chosen active sets---especially during the first $10$ MPC iterations. However, for $\overline{m} \geq 10$ a reasonably accurate approximation of the optimal number of active constraints is maintained during all iterations.}

\begin{figure}[htbp!]
    \centering
    \includegraphics[width=0.41\textwidth]{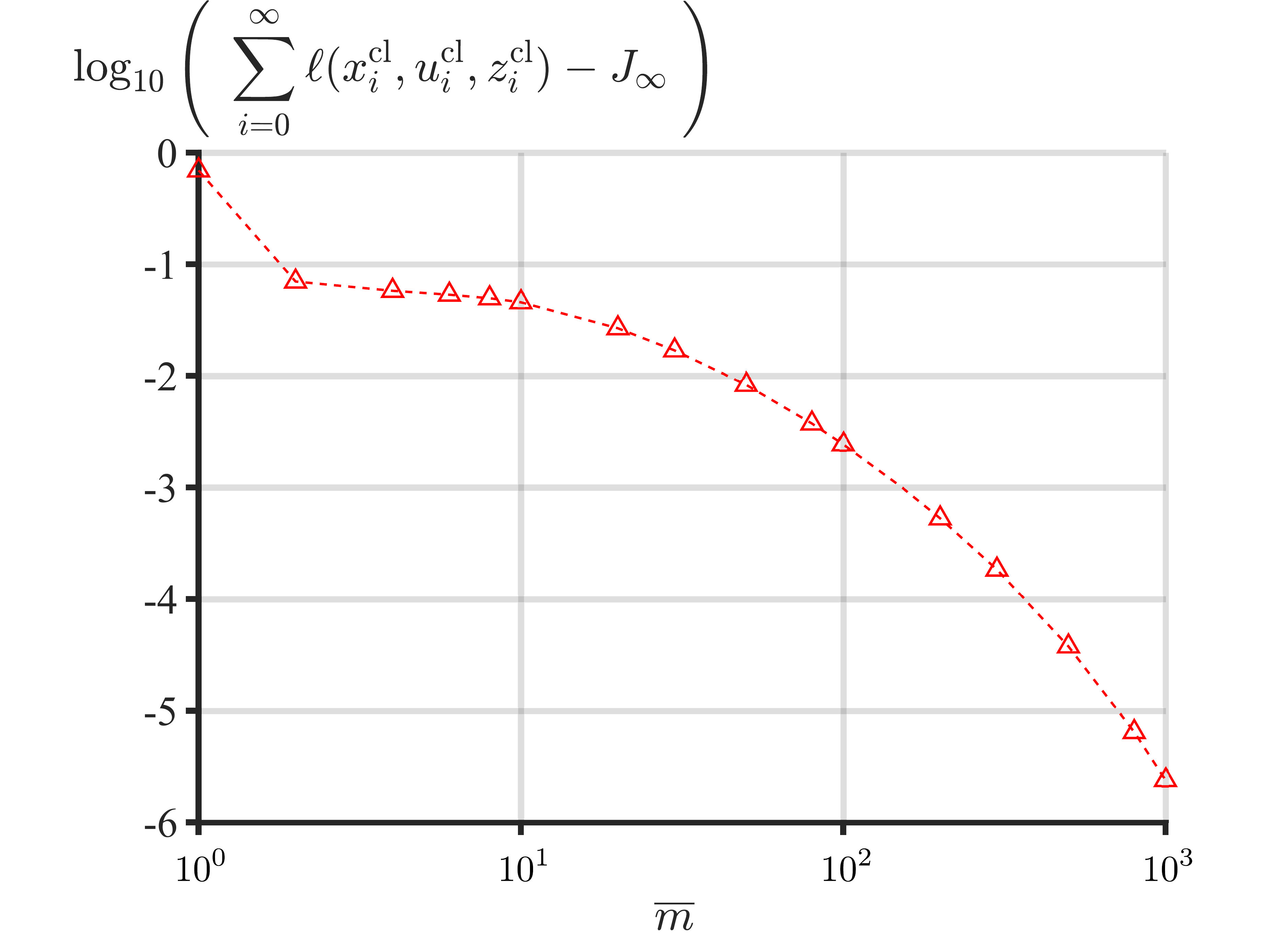}
    \caption{Closed-loop performance degradation (log scale) with
        respect to the optimal objective function $J_{\infty}$ as a
        function of the number of iterations $\overline{m}$ in Algorithm~\ref{alg::dempc}.}
    \label{fig::close} 
\end{figure}

Finally, Figure~\ref{fig::close} shows the sub-optimality of Algorithm~1 in dependence on $\overline{m}$ for a representative case study with $\bar I = 3$ and $N=30$---for other values of $\bar I$ and $N$ the convergence looks similar.

\section{Conclusions}
\label{sec::conclusions}
This paper has introduced a parallelizable and real-time verifiable MPC scheme, presented in the form of Algorithm~\ref{alg::dempc}. This control algorithm evaluates at every sampling time a finite number of pre-computed, explicit piecewise affine solution maps that are associated with parametric small-scale QPs. Because solving large-scale QPs in real-time may be impossible, the presented algorithm returns suboptimal control reaction on purpose---in order to be able to meet challenging real-time and limited memory requirements. The theoretical contributions of this paper have been presented in Theorem~\ref{thm::stability} and Corollary~\ref{cor::performance}, which provide both \change{asymptotic stability} guarantees as well as bounds on sub-optimality. The presented explicit MPC approach can be used to reduce the storage and run-time of explicit MPC by orders of magnitude, as illustrated by applying Algorithm~1 to a spring-mass-vehicle benchmark problem.

\begin{ack}
    \tiny
    \change{The authors thank three anonymous reviewers and the associate editor for their precise and substantial remarks, which have helped to significantly improve this article.}
    All authors were supported via a bilateral grant agreement between
    China and Slovakia (SK-CN-2015-PROJECT-6558, APVV SK-CN-2015-0016), 
    which is gratefully acknowledged. The work of Yuning Jiang and Boris Houska
    was supported by the National Science Foundation China (NSFC), 
    Nr. 61473185, as well as ShanghaiTech University, Grant-Nr. 
    F-0203-14-012. Juraj Oravec and Michal Kvasnica acknowledge the 
    contribution of the Scientific Grant Agency of the Slovak 
    Republic under the grants 1/0585/19, 1/0112/16, and the Slovak Research 
    and Development Agency under the project APVV-15-0007.
\end{ack}


\begin{thebibliography}{10}
	
	\bibitem{BemEtal:aut:02}
	A.~Bemporad, M.~Morari, V.~Dua, and E.N. Pistikopoulos.
	\newblock The explicit linear quadratic regulator for constrained systems.
	\newblock {\em Automatica}, 38(1):3--20, 2002.
	
	\bibitem{Bertsekas2012}
	D.P. Bertsekas.
	\newblock {\em Dynamic Programming and Optimal Control}.
	\newblock Athena Scientific Dynamic Programming and Optimal Control, Belmont,
	Massachusetts, 3rd edition, 2012.
	
	\bibitem{Bla:aut:99}
	F.~Blanchini.
	\newblock {Set invarinace in control --- A survey}.
	\newblock {\em Automatica}, 35:1747--1767, 1999.
	
	\bibitem{BorrelliPHD}
	F.~Borrelli.
	\newblock {\em Constrained Optimal Control Of Linear And Hybrid Systems},
	volume 290 of {\em Lecture Notes in Control and Information Sciences}.
	\newblock Springer, 2003.
	
	\bibitem{borrelli2003geometric}
	F.~Borrelli, A.~Bemporad, and M.~Morari.
	\newblock Geometric algorithm for multiparametric linear programming.
	\newblock {\em Journal of optimization theory and applications},
	118(3):515--540, 2003.
	
	\bibitem{Borrelli2017}
	F.~Borrelli, A.~Bemporad, and M.~Morari.
	\newblock {\em Predictive control for linear and hybrid systems}.
	\newblock Cambridge University Press, 2017.
	
	\bibitem{Boyd2011}
	S.~Boyd, N.~Parikh, E.~Chu, B.~Peleato, and J.~Eckstein.
	\newblock Distributed optimization and statistical learning via the alternating
	direction method of multipliers.
	\newblock {\em Foundation Trends in Machine Learning}, 3(1):1--122, 2011.
	
	\bibitem{BoyVan:ConOpt:04}
	S.~Boyd and L.~Vandenberghe.
	\newblock {\em Convex Optimization}.
	\newblock Cambridge University Press, 2004.
	
	\bibitem{Cagienard2007}
	R.~Cagienard, P.~Grieder, E.C. Kerrigan, and M.~Morari.
	\newblock Move blocking strategies in receding horizon control.
	\newblock {\em Journal of Process Control}, 17(6):563--570, 2007.
	
	\bibitem{Conte2016}
	C.~Conte, C.N. Jones, M.~Morari, and M.N. Zeilinger.
	\newblock Distributed synthesis and stability of cooperative distributed model
	predictive control for linear systems.
	\newblock {\em Automatica}, 69:117--125, 2016.
	
	\bibitem{Conte2012}
	C.~Conte, T.~Summers, M.N. Zeilinger, M.~Morari, and C.N. Jones.
	\newblock Computational aspects of distributed optimization in model predictive
	control.
	\newblock In {\em Proceedings of the 51st IEEE Conference on Decision and
		Control, 2012}, pages 6819--6824, 2012.
	
	\bibitem{Diehl2002}
	M.~Diehl, H.G. Bock, J.P. Schl{\"o}der, R.~Findeisen, Z.~Nagy, and
	F.~Allg{\"o}wer.
	\newblock Real-time optimization and nonlinear model predictive control of
	processes governed by differential-algebraic equations.
	\newblock {\em Journal of Process Control}, 12(4):577--585, 2002.
	
	\bibitem{Diehl2009}
	M.~Diehl, H.J. Ferreau, and N.~Haverbeke.
	\newblock Efficient numerical methods for nonlinear {MPC} and moving horizon
	estimation.
	\newblock In L.~Magni, M.D. Raimondo, and F.~Allg\"ower, editors, {\em
		Nonlinear model predictive control}, volume 384 of {\em Lecture Notes in
		Control and Information Sciences}, pages 391--417. Springer, 2009.
	
	\bibitem{Everett1963}
	H.~Everett.
	\newblock Generalized {L}agrange multiplier method for solving problems of
	optimum allocation of resources.
	\newblock {\em Operations Research}, 11(3):399--417, 1963.
	
	\bibitem{FS16}
	L.~Ferranti, G.~Stathopoulos, C.~N. Jones, and T.~Keviczky.
	\newblock Constrained {LQR} using online decomposition techniques.
	\newblock In {\em IEEE Conference on Decision and Control}, pages 2339--2344,
	Las Vegas, USA, 2016.
	
	\bibitem{Ferreau2014}
	H.~J. Ferreau, C.~Kirches, A.~Potschka, H.~G. Bock, and M.~Diehl.
	\newblock qpoases: A parametric active-set algorithm for quadratic programming.
	\newblock {\em Mathematical Programming Computation}, 6(4):327--363, 2014.
	
	\bibitem{Giselesson2013}
	P.~Giselsson, M.~Dang~Doan, T.~Keviczky, B.~De~Schutter, and A.~Rantzer.
	\newblock Accelerated gradient methods and dual decomposition in distributed
	model predictive control.
	\newblock {\em Automatica}, 49(3):829--833, 2013.
	
	\bibitem{GM02}
	P.~Grieder and M.~Morari.
	\newblock Complexity reduction of receding horizon control.
	\newblock In {\em IEEE Conference on Decision and Control}, Maui, Hawaii, USA,
	2003.
	
	\bibitem{Gruene2009}
	L.~Gr\"une.
	\newblock Analysis and design of unconstrained nonlinear mpc schemes for finite
	and infinite dimensional systems.
	\newblock {\em SIAM Journal on Control and Optimization}, 48(2):1206--1228,
	2009.
	
	\bibitem{HJ15}
	M.~Herceg, C.~N. Jones, M.~Kvasnica, and M.~Morari.
	\newblock Enumeration-based approach to solving parametric linear
	complementarity problems.
	\newblock {\em Automatica}, (62):243--248, 2015.
	
	\bibitem{MPT3}
	M.~Herceg, M.~Kvasnica, C.~Jones, and M.~Morari.
	\newblock Multi-parametric toolbox 3.0.
	\newblock In {\em 2013 European Control Conference}, pages 502--510, 2013.
	
	\bibitem{Houska2011}
	B.~Houska, H.J. Ferreau, and M.~Diehl.
	\newblock An auto-generated real-time iteration algorithm for nonlinear {MPC}
	in the microsecond range.
	\newblock {\em Automatica}, 47:2279–2285, 2011.
	
	\bibitem{Houska2016}
	B.~Houska, J.~Frasch, and M.~Diehl.
	\newblock An augmented {L}agrangian based algorithm for distributed non-convex
	optimization.
	\newblock {\em SIAM Journal on Optimization}, 26(2):1101--1127, 2016.
	
	\bibitem{Ingole2015}
	D.~Ingole and M.~Kvasnica.
	\newblock {FPGA} implementation of explicit model predictive control for closed
	loop control of depth of anesthesia.
	\newblock In {\em 5th IFAC Conference on Nonlinear Model Predictive Control},
	pages 484--489, 2015.
	
	\bibitem{KL18}
	B.~Karg and S.~Lucia.
	\newblock Efficient representation and approximation of model predictive
	control laws via deep learning.
	\newblock {\em arXiv.org/math.OC}, abs/1806.10644, 2018.
	
	\bibitem{Kozma2013}
	A.~Kozma, J.V. Frasch, and M.~Diehl.
	\newblock A distributed method for convex quadratic programming problems
	arising in optimal control of distributed systems.
	\newblock In {\em Proceedings of the IEEE Conference on Decision and Control
		(CDC)}, 2013.
	
	\bibitem{Kva:regionless:2015}
	M.~Kvasnica, B.~Tak\'acs, J.~Holaza, and S.~Di~Cairano.
	\newblock On region-free explicit model predictive control.
	\newblock In {\em Conference on Decision and Control (CDC)}, pages 3669--3674,
	2015.
	
	\bibitem{yalmip}
	J.~L{\"o}fberg.
	\newblock {YALMIP}, 2004.
	\newblock Available from \url{http://users.isy.liu.se/johanl/yalmip/}.
	
	\bibitem{Mattingley2009}
	J.~Mattingley and S.~Boyd.
	\newblock Automatic code generation for real-time convex optimization.
	\newblock {\em Convex optimization in signal processing and communications},
	pages 1--41, 2009.
	
	\bibitem{Necoara2008}
	I.~Necoara and J.A.K. Suykens.
	\newblock Application of a smoothing technique to decomposition in convex
	optimization.
	\newblock {\em IEEE Transactions on Automatic Control}, 53(11):2674--2679,
	2008.
	
	\bibitem{oberdieck2016multi}
	R.~Oberdieck, N.~Diangelakis, I.~Nascu, M.~Papathanasiou, M.~Sun,
	S.~Avraamidou, and E.~Pistikopoulos.
	\newblock On multi-parametric programming and its applications in process
	systems engineering.
	\newblock {\em Chemical Engineering Research and Design}, 116:61--82, 2016.
	
	\bibitem{Donoghue2013}
	B.~O'Donoghue, G.~Stathopoulos, and S.~Boyd.
	\newblock A splitting method for optimal control.
	\newblock {\em IEEE Transactions on Control Systems Technology},
	21(6):2432--2442, 2013.
	
	\bibitem{Oravec2017}
	J.~Oravec, Y.~Jiang, B.~Houska, and M.~Kvasnica.
	\newblock Parallel explicit {MPC} for hardware with limited memory.
	\newblock In {\em In Proceedings of the 20th IFAC World Congress, Toulouse,
		France}, pages 3356--3361, 2017.
	
	\bibitem{Qin2003}
	S.J. Qin and T.A. Badgwell.
	\newblock A survey of industrial model predictive control technology.
	\newblock {\em Control Engineering Practice}, 93(316):733--764, 2003.
	
	\bibitem{Rawlings2009}
	J.B. Rawlings, D.Q. Mayne, and M.M. Diehl.
	\newblock {\em Model Predictive Control: Theory and Design}.
	\newblock Madison, WI: Nob Hill Publishing, 2018.
	
	\bibitem{Richter2011}
	S.~Richter, M.~Morari, and C.N. Jones.
	\newblock Towards computational complexity certification for constrained {MPC}
	based on {L}agrange relaxation and the fast gradient method.
	\newblock In {\em Proceedings of the 50th IEEE Conference on Decision and
		Control and European Control Conference}, pages 5223--5229, 2011.
	
	\bibitem{SK13}
	G.~Stathopoulos, T.~Keviczky, and Y.~Wang.
	\newblock A hierarchical time-splitting approach for solving finite-time
	optimal control problems.
	\newblock In {\em European Control Conference}, pages 3089--3094, Z\"{u}rich,
	Switzerland, 2013.
	
	\bibitem{TJB03}
	P.~T{\o}ndel, T.A. Johansen, and A.~Bemporad.
	\newblock Evaluation of piecewise affine control via binary search tree.
	\newblock {\em Automatica}, 39(5):945--950, 2003.
	
	\bibitem{Wang2015}
	Y.~Wang, B.~O'Donoghue, and S.~Boyd.
	\newblock Approximate dynamic programming via iterated {B}ellman inequalities.
	\newblock {\em International Journal of Robust and Nonlinear Control},
	25(10):1472--1496, 2015.
	
	\bibitem{Zavala2009}
	V.M. Zavala and L.T. Biegler.
	\newblock The advanced-step {NMPC} controller: Optimality, stability and
	robustness.
	\newblock {\em Automatica}, 45(1):86--93, 2009.
	
\end{thebibliography}

\appendix

\footnotesize

\section{Proof of Theorem~\ref{thm::convergenceRate}}
\label{app::convergenceRate}
\change{As the proof of Theorem~\ref{thm::convergenceRate} is not entirely straightforward,
we first establish an intermediate results on the convergence properties
of Algorithm~\ref{alg::dempc}, which is then, in a second step, used to obtain
a convergence rate estimate. Therefore, this proof is divided into two subsections:
Subsection~\ref{app::convergenceProperties} analyzes the general convergence
properties of Algorithm~\ref{alg::dempc} and Subsection~\ref{app::convergenceRate1}
presents a complete proof of Theorem~\ref{thm::convergenceRate} by using
these properties.}

\subsection{\change{A closer look at the convergence properties of Algorithm~\ref{alg::dempc}}}
\label{app::convergenceProperties}

\change{The goal of the section is to establish the following technical result.}
\begin{lemma}
\label{lem::convergence}
Let Assumption~\ref{ass::blanket} be satisfied and let~\eqref{eq::mpc2} be feasible, i.e., such that a unique minimizer $y\opt$ and an associated dual solution $\lambda\opt$ exist. Then the iterates of Algorithm~\ref{alg::dempc} satisfy
\[
\sum_{m=\hat m}^{\overline{m}} F( \xi^m - y\opt ) \leq \frac{ F( y^{\hat m} - y\opt ) + F\opt( \lambda^{\hat m} - \lambda\opt ) }{4}
\]
for all $\overline{m} \geq \hat m$ and all $\hat m \geq 2$.
\end{lemma}

\change{\textbf{Proof.}} Let us introduce the auxiliary functions
\begin{equation}
\begin{split}\notag
\mathcal F_0(\phi_0) \;=\;& F_0(\phi_0)  - \left( H_0^\tr \lambda_{0}^m \right)^\tr \phi_0^m + \nabla F_0( \xi_0^m - y_0^m )^\tr \phi_0 \,, \\[0.16cm]
\mathcal F_k(\phi_k) \;=\;& F_k(\phi_k)  + \left( G_k^\tr \lambda_{k-1}^{m} - H_k^\tr \lambda_{k}^m \right)^\tr \phi_k^m \\
& + \nabla F_k( \xi_k^m - y_k^m )^\tr \phi_k \; . 
\end{split}
\end{equation}
Because $\xi_k^m$ is a minimizer of the $k$-th decoupled QP in Step 2a) of Algorithm~1, it must also be a minimizer of $\mathcal F_k$ on $\mathbb Y_k$. Thus, because $\mathcal F_k$ is strongly convex with Hessian $\nabla^2 F_k$, we must have
\[
\sum_{k=0}^N \mathcal F_k(\xi_k^m) + \sum_{k=0}^N F_k(\xi_k^m - y_k\opt) \leq \sum_{k=0}^N \mathcal F_k(y_k\opt)\,.
\]
On the other hand, due to duality, we have
\begin{align}
&\sum_{k=0}^N F_k(y_k\opt) + \langle \lambda\opt , y\opt \rangle + \sum_{k=0}^N F_k(\xi_k^m - y_k\opt) \notag \\
\leq\;\; &\sum_{k=0}^N F_k(\xi_k^m) + \langle \lambda\opt , \xi^m \rangle , \notag
\end{align}
where the shorthand notation
\[
\langle \lambda , y \rangle = - \left( H_0^\tr \lambda_{0} \right)^\tr y_0 + \sum_{k=1}^N \left( G_k^\tr \lambda_{k-1} - H_k^\tr \lambda_{k} \right)^\tr y_k
\]
is used to denote a weighted (non-symmetric) scalar product of primal and dual variables. Adding both inequalities and collecting terms yields
\begin{equation}\notag
\begin{split}
0 \;\geq\;& \sum_{k=0}^N \nabla F_k( \xi_k^m - y_k^m )^\tr ( \xi_k^m - y_k\opt ) + 2 \sum_{k=0}^N F_k(\xi_k^m - y_k\opt) \\[0.16cm]
& + \langle \lambda^m - \lambda\opt, \xi^m - y\opt \rangle \; . 
\end{split}
\end{equation}
Let us introduce the matrices 
\[
\mathcal Q = \nabla^2  \left( \sum_{k=0}^N  F_k\right)  \;,\;\mathcal A = \nabla_{\lambda,x} \langle  \lambda,y \rangle
\] 
such that the above inequality can be written in the form
\begin{align}
0\; \geq\;&\;\; (\xi^m - y^m )^\tr \mathcal Q ( \xi^m - y\opt ) + 2 \sum_{k=0}^N F_k(\xi_k^m - y_k\opt) \notag \\[0.16cm]
& \;\;+ \left( \lambda^m - \lambda\opt \right)^\tr \mathcal A \left( \xi^m - y\opt \right) \; .
\label{eq::aux1}
\end{align}
Similarly, the stationarity condition QP~\eqref{eq::qp} can be written as
\[
\mathcal Q( y^{m+1} - 2 \xi^m + y^m ) + \mathcal A^\tr \delta^m = 0 \; .
\]
Because $\mathcal Q$ is positive definite, we solve this equation with respect to $\xi^m$ finding
\begin{equation}
\label{eq::QPopt}
\xi^m = \frac{1}{2} \mathcal Q^{-1} \mathcal A^\tr (\lambda^{m+1} - \lambda^m) + \frac{y^m + y^{m+1}}{2}\;.
\end{equation}
Here, we have additionally substituted the relation
$$\delta^m = \lambda^{m+1} - \lambda^m \; .$$
Notice that we have $\mathcal A y^m = \mathcal A y^{m+1} = \mathcal A y\opt$ for all $m \geq 2$, because the solutions of the QP~\eqref{eq::qp} must satisfy the equality constraints in~\eqref{eq::mpc2}. If we substitute this equation and the expression for $\xi^m$ in~\eqref{eq::aux1}, we find that
\begin{align}\notag
&\;- 2 F(\xi^m - y\opt) \\\notag
\geq &\;\; (\xi^m - y^m )^\tr \mathcal Q ( \xi^m - y\opt ) + \left( \lambda^m - \lambda\opt \right)^\tr \mathcal A \left( \xi^m - y\opt \right) \\[0.16cm]\notag
=&\;\; \frac{1}{4}(\lambda^{m+1} - \lambda^m)^\tr \mathcal A \mathcal Q^{-1} \mathcal A^\tr (\lambda^{m+1} - \lambda^m) \\[0.16cm]
\label{eq::descent}
& \;\;+ \frac{1}{4}(y^{m+1} - y^{m}) \mathcal Q ( y^m - 2y\opt + y^{m+1} ) \\[0.16cm]\notag
&\;\; + \frac{1}{2}(\lambda^{m} - \lambda\opt)^\tr \mathcal A \mathcal Q^{-1} \mathcal A^\tr (\lambda^{m+1} - \lambda^m)\\[0.16cm]\notag
=& \;\;\frac{1}{2} \left( F(y^{m+1} - y\opt )  - F(y^{m} - y\opt ) \right) \\[0.16cm]\notag
& \;\;+  \frac{1}{2} \left( F\opt( \lambda^{m+1} - \lambda\opt ) - F\opt( \lambda^{m} - \lambda\opt ) \right)
\end{align}
for all $m \geq 2$. Now, the statement of \change{Lemma~\ref{lem::convergence}} follows by summing up the above inequalities for $m=\hat m$ to $m=\overline{m}$ and using that the last element in the telescoping sum on the right hand,
\[
\frac{F(y^{\overline{m}+1} - y\opt ) + F\opt( \lambda^{\overline{m}+1} - \lambda\opt )}{2} \geq 0 
\]
is non-negative.

\subsection{Analysis of the convergence rate of Algorithm~\ref{alg::dempc}}
\label{app::convergenceRate1}

\change{The goal of this section is to prove the statement of Theorem~\ref{thm::convergenceRate} by using the intermediate result from the previous section.}
Let $\hat{\mathbb Y}_k$ denote the intersection of all active supporting hyperplanes at the solutions of the small scale QPs of Step 2a) in Algorithm~\ref{alg::dempc} for $k \in \{ 0, \ldots, N\}$ at a given iteration $m$. We construct the auxiliary optimization problem\\[-0.45cm]
\begin{align}
\label{eq::mpc2Aux}\notag
\underset{\hat y}{\min}&\;\; \sum_{k=0}^{N} F_k( \hat y_k ) \\[0.16cm]
\text{s.t.} &
\left\{
\begin{array}{l}
\forall k \in \{ 0, \ldots, N-1 \}, \\
\begin{array}{rcll}
G_{k+1} \hat y_{k+1} &=& H_k \hat y_k + h_k \; & \mid \; \hat \lambda_{k} \; , \\
0 &=& H_N \hat y_N \; & \mid \; \hat \lambda_{N} \; ,
\end{array} \\
\hat y_k \in \hat{\mathbb Y}_k \; , \; \hat y_N \in \hat{\mathbb Y}_N \; ,
\end{array}
\right.
\end{align}
and denote optimal primal and dual solutions of this problem by $\hat y\opt$ and $\hat \lambda\opt$. Next, we also construct the auxiliary QPs
\begin{align}
\underset{\xi_0^m \in \hat{\mathbb Y}_0}{\min} & \; F_0(\xi_0^m)  - \left( H_0^\tr \lambda_{0}^m \right)^\tr \xi_0^m + F_0( \xi_0^m - y_0^m ) , \notag \\[0.2cm]
\underset{\xi_k^m \in \hat{\mathbb Y}_k}{\min} & \; F_k(\xi_k^m)  + \left( G_k^\tr \lambda_{k-1}^m - H_k^\tr \lambda_{k}^m \right)^\tr \xi_k^m + F_k( \xi_k^m - y_k^m )  \; . \notag
\end{align}
Because these QPs have equality constraints only, their parametric solutions must be affine. Thus, there exists a matrix $T_1$ such that
\[
\xi^m - \hat y\opt = T_1 \left(
\begin{array}{c}
y^m - \hat y\opt \\
\lambda^m - \hat \lambda\opt
\end{array}
\right).
\]
Similarly, the coupled QP~\eqref{eq::qp} has equality constraints only, i.e., there exists a matrix $T_2$ such that
\[
\left(
\begin{array}{c}
y^{m+1} - \hat y\opt \\
\delta^m \\
\end{array}
\right) = T_2 \left(
\begin{array}{c}
\xi^m - \hat y\opt \\
y^m - \hat y\opt
\end{array}
\right).
\]
Now, we use the equation $\lambda^{m+1} - \lambda\opt = \lambda^m - \lambda\opt + \delta$ and substitute the above equations finding that
\begin{equation}
\left(
\begin{array}{c}
y^{m+1} - \hat y\opt \\
\lambda^{m+1} - \hat \lambda\opt
\end{array}
\right) = T
\left(
\begin{array}{c}
y^{m} - \hat y\opt \\
\lambda^{m} - \hat \lambda\opt
\end{array}
\right)
\label{eq::linearOperator}
\end{equation}
with
\[
\quad T = \left(
\begin{array}{c}
T_2 \left(
\begin{array}{c}
T_1 \\
(I \;\; 0 )
\end{array}
\right) + (0 \;\; I)
\end{array}
\right) \; .
\]
Next, we know from \change{Lemma~\ref{lem::convergence}} that if we would apply Algorithm~\ref{alg::dempc} to the auxiliary problem~\eqref{eq::mpc2Aux}, the corresponding primal and dual iterates would converge to $\hat y\opt$ and $\hat \lambda\opt$. In particular, inequality~\eqref{eq::descent} from the proof of \change{Lemma~\ref{lem::convergence}} can be applied finding that
\begin{equation}
\begin{array}{cl}
&\left( y^{m+1} - \hat y\opt \right)^\tr \mathcal Q \left( y^{m+1} - \hat y\opt \right) \\
& + \left( \lambda^{m+1} - \hat \lambda\opt \right)^\tr  \mathcal{A} \mathcal{Q}^{-1} \mathcal{A}^\tr \left( \lambda^{m+1} - \hat \lambda\opt \right) \\
<& \left( y^{m} - \hat y\opt \right)^\tr \mathcal Q \left( y^{m} - \hat y\opt \right) \\
& + \left( \lambda^{m} - \hat \lambda\opt \right)^\tr \mathcal A Q^{-1} \mathcal A^\tr \left( \lambda^{m} - \hat \lambda\opt \right) \; ,
\end{array}
\end{equation}
whenever $\left(
\begin{array}{c}
y^{m} - \hat y\opt\\
\lambda^{m} - \hat \lambda\opt 
\end{array} 
\right) 
\neq 0$. By substituting the linear equation~\eqref{eq::linearOperator}, we find that this is only possible if
\[
T^\tr \left(
\begin{array}{cc}
\mathcal Q & 0 \\
0 & \mathcal A Q^{-1} \mathcal A^\tr
\end{array}
\right)
T \preceq \kappa_{\mathbb A} I
\]
for a constant $\kappa_{\mathbb A} < 1$. Now, one remaining difficulty is that the constant $\kappa_{\mathbb A}$ (as well as the matrix $T$) depends on the particular set $\mathbb A$ of active supporting hyperplanes in the small-scale QPs. Nevertheless, because there exists only a finite number of possible active sets, the maximum
\[
\kappa = \max_{\mathbb A} \kappa_{\mathbb A}
\]
must exist and satisfy $\kappa < 1$. Now, the equation
\begin{equation}
\left(
\begin{array}{c}
y^{m+1} - y\opt \\
\lambda^{m+1} - \lambda\opt
\end{array}
\right) = T
\left(
\begin{array}{c}
y^{m} - y\opt \\
\lambda^{m} - \lambda\opt
\end{array}
\right) 
\end{equation}
holds only for our fixed $m$ and the associated matrix $T$ for a particular active set, but the associated decent condition
\begin{equation}
\begin{array}{l}
\quad\left( y^{m+1} - y\opt \right)^\tr \mathcal Q \left( y^{m+1} - y\opt \right) \\
\quad + \left( \lambda^{m+1} - \lambda\opt \right)^\tr \mathcal A Q^{-1} \mathcal A^\tr \left( \lambda^{m+1} - \lambda\opt \right) \\
\leq \kappa \left[ \left( y^{m} - y\opt \right)^\tr \mathcal Q \left( y^{m} - y\opt \right) \right. \\
\left. \quad + \left( \lambda^{m} - \lambda\opt \right)^\tr \mathcal A Q^{-1} \mathcal A^\tr \left( \lambda^{m} - \lambda\opt \right) \right] \; ,
\end{array}
\end{equation}
holds independently of the active set of the QPs in the $m$-th iteration and is indeed valid for all $m$. After re-introducing the functions $F$ and $F\opt$, we obtain the statement of the theorem.

\StartChange

\section{Proof of Proposition~\ref{prop::feasibility}}
\label{app::feasibility}
Because $\xi_0^{\overline{m}}$ is a feasible solution of the first small-scale decoupled QP,
we have $\xi_0^{\overline{m}} \in \mathbb Y_0$. Notice that such a feasible solution exists
due to the particular construction of the set $\mathbb Y_0$ and our assumption that
the set $\mathbb X$ is control invariant. Consequently, the next iterate for the state,
\[
x_0^+ = A x_0 + (B \;\; C) \xi_0^{\overline{m}} 
\]
must satisfy $x_0^+ \in \mathbb X$ by construction.
This is the first statement of Proposition~\ref{prop::feasibility}.
In order to establish the second statement, we observe that the particular construction
of the sets $\mathbb Y_k$ implies that there exists a feasible point of~\eqref{eq::mpc2} for
any choice of $x_0 \in \mathbb X$, because we choose $\mathbb X_N = \mathbb X$, but
$\mathbb X$ is control invariant. Thus, we must have
$\mathbb X \subseteq \mathcal X$. The other way around, if $x_0 \notin \mathbb X$
the state constraints are violated, i.e., we also have $\mathcal X \subseteq \mathbb X$.
Thus, we have $\mathcal X = \mathbb X$ and, consequently, $\mathcal X$
is a control invariant set, too. This is sufficient to establish recursive feasibility of Algorithm~\ref{alg::dempc}.
\EndChange

\section{Proof of Theorem~\ref{thm::stability}}
\label{app::stability}

\change{Recall that stability proofs for standard MPC proceed by using the inequality
\begin{equation}
\label{eq::MPC_Lyaponov}
J(x_1\opt) \leq J(x_0) - F_0( y_0\opt ) \; ,
\end{equation}
which holds if Assumption~\ref{ass::TerminalLyapunov} is satisfied (see, e.g.,~\cite{Gruene2009} for a derivation of this inequality), i.e., $J$ is a global Lyapunov function on $\mathcal X$. Now, if we implement Algorithm~\ref{alg::dempc} with a finite $\overline{m}$, we have
\begin{equation}
\label{eq::lyapD1}
J(x_0^+) \leq J(x_0) - \left( F_0( y_0\opt ) - J(x_0^+) + J(x_1\opt)  \right) \; ,
\end{equation}
i.e., $J$ can still be used as a Lyaponov function proving asymptotic stability as long as we ensure that
\begin{equation}
\label{eq::lyapD2}
F_0( y_0\opt ) - J(x_0^+) + J(x_1\opt) \geq \alpha F_0( y_0\opt ) \; ,
\end{equation}
for a constant $\alpha > 0$. In order to show that such an inequality can indeed be established for a finite $\overline{m}$, we need the following technical result.}

\begin{lemma}
\label{lem::bound}
\change{The iterate $x_0^+$ satisfies the inequality}
\[
\left( x_0^+ - x_1\opt \right)^\tr Q \left( x_0^+ - x_1\opt \right) \leq \sigma \gamma^2 ( 1 + \sqrt{\kappa})^2 \kappa^{\overline{m}} F_0( y_0\opt ) \; .
\]
\end{lemma}

\textbf{Proof.} We start with the equation
\[
x_0^+ - x_1\opt = (B\;\;C) ( \xi_1^{\overline{m}} - y_0\opt ) = \mathcal P \mathcal A \left( \xi^{\overline{m}} - y\opt \right) \; ,
\]
which holds for the projection matrix $\mathcal P = \mathrm{diag}(I , 0 ,$ $\ldots, 0 )$ that filters out the first block component of the equality constraint residuum $\mathcal A \left( \xi^{\overline{m}} - y\opt \right)$. Next, we substitute~\eqref{eq::QPopt}, which yields
\begin{align}
& \;x_0^+ - x_1\opt \notag \\ 
=& \;\mathcal P \mathcal A \left[ \frac{\mathcal Q^{-1} \mathcal A^\tr (\lambda^{{\overline{m}}+1} - \lambda^{\overline{m}})}{2} + \frac{y^{\overline{m}+1} + y^{{\overline{m}}}}{2} - y\opt \right] \notag \\
=&\; \frac{1}{2} \mathcal P \mathcal A \mathcal Q^{-1} \mathcal A^\tr (\lambda^{{\overline{m}}+1} - \lambda^{\overline{m}}) \; ,
\end{align}
where we have used that $0 = \mathcal A( y^{\overline{m}} - y\opt )$ and $0 = \mathcal A( y^{{\overline{m}}+1} - y\opt )$ recalling that these equations follow from the equality constraints in the coupled QP in Step 2b) of Algorithm~\ref{alg::dempc}. The particular definition of $\sigma$ implies that
\begin{align}
&4 \left( x_0^+ - x_1\opt \right)^\tr Q \left( x_0^+ - x_1\opt \right) \notag \\[0.16cm]
\leq& (\lambda^{{\overline{m}}+1} - \lambda^{\overline{m}})^\tr 
\mathcal A \mathcal Q^{-1} \underbrace{\mathcal A^\tr \mathcal P^\tr Q \mathcal P \mathcal A}_{\preceq \sigma \mathcal Q} \mathcal Q^{-1} \mathcal A^\tr (\lambda^{{\overline{m}}+1} - \lambda^{\overline{m}}) \notag \\[0.16cm]
\leq& \sigma (\lambda^{{\overline{m}}+1} - \lambda^{\overline{m}})^\tr \mathcal A \mathcal Q^{-1} \mathcal A^\tr (\lambda^{{\overline{m}}+1} - \lambda^{\overline{m}}) \; . \notag
\end{align}
Now, we can use the result of Theorem~\ref{thm::convergenceRate} to find
\begin{align}
&\quad (\lambda^{{\overline{m}}+1} - \lambda^{\overline{m}})^\tr \mathcal A \mathcal Q^{-1} \mathcal A^\tr (\lambda^{{\overline{m}}+1} - \lambda^{\overline{m}}) \notag \\
&\leq F\opt(\lambda^{{\overline{m}}+1}-\lambda\opt) + F\opt(\lambda^{{\overline{m}}}-\lambda\opt) \notag \\
& \quad + 2 \sqrt{F\opt(\lambda^{{\overline{m}}+1}-\lambda\opt) F\opt(\lambda^{{\overline{m}}}-\lambda\opt) } \notag \\
&\leq \kappa^{\overline{m}} ( 1 + \sqrt{\kappa})^2 \left( F(y^{1}-y\opt) + F\opt(\lambda^{1}-\lambda\opt) \right)  \; . \notag
\end{align}
It remains to use the inequalities
\begin{align}
F(y\opt) + F\opt(\lambda\opt) &\leq \gamma^2 x_0^\tr Q x_0 \leq \gamma^2 F_0( y_0\opt )\;, \notag \\
F(y^1) + F\opt(\lambda^1) &\leq \gamma^2 x_0^\tr Q x_0 \leq \gamma^2 F_0( y_0\opt ) \; , \notag
\end{align}
which hold due to Assumption~\ref{ass::gamma} and the particular construction in Step~1 of Algorithm~\ref{alg::dempc}, arriving at the inequality
\begin{align}\notag
&\;\;4 \left( x_0^+ - x_1\opt \right)^\tr Q \left( x_0^+ - x_1\opt \right) \\\notag
\leq &\;\;\sigma (\lambda^{{\overline{m}}+1} - \lambda^{\overline{m}})^\tr \mathcal A \mathcal Q^{-1} \mathcal A^\tr (\lambda^{{\overline{m}}+1} - \lambda^{\overline{m}}) \\\notag
\leq &\;\;\sigma \kappa^{\overline{m}} ( 1 + \sqrt{\kappa})^2 \left( F(y^{1}-y\opt) + F\opt(\lambda^{1}-\lambda\opt) \right) \\\notag
\leq&\;\; 4 \sigma \kappa^{\overline{m}} ( 1 + \sqrt{\kappa})^2 \gamma^2 F_0{y_0\opt} \; .
\end{align}
Dividing by $4$ on both sides yields the statement of the lemma. \hfill\qed

\change{By combining the inequality from the above lemma and inequality~\eqref{eq::etaTauBound} we find
\begin{align}
\label{eq::auxA}
&\left|J(x_0^+) - J(x_1\opt)\right|  \\\notag
\leq& \left[ \eta \sqrt{\sigma} \gamma (1+\sqrt{\kappa}) + \frac{\tau \sigma \gamma^2 (1+\sqrt{\kappa})^2}{2} \right] \kappa^{\frac{\overline{m}}{2}} F_0( y_0\opt ) \; .
\end{align}
Now, the statement of Theorem~\ref{thm::stability} follows directly from~\eqref{eq::auxA}, as the construction of $\overline{m}$ ensures that the Lyapunov descent condition~\eqref{eq::lyapD2} holds on $\mathcal X$ with
\[
\alpha = 1 - \left[ \eta \sqrt{\sigma} \gamma (1+\sqrt{\kappa}) + \frac{\tau \sigma \gamma^2 (1+\sqrt{\kappa})^2}{2} \right] \kappa^{\frac{\overline{m}}{2}} > 0 \; ,
\]
which is sufficient to establish asymptotic stability.}

\end{document}